\documentstyle[11pt]{article}
\begin{document}
\begin{center}
{\Large\bf CAN YOU FEEL THE DOUBLE JUMP \\ Saharon Shelah \\ Joel Spencer
\\ {\today} \\ Preliminary Draft}
\end{center}
\vspace{2 cm}
\section{Summary of Results}
 In their fundamental work Paul Erd\H{o}s and Alfred Renyi [1] 
considered the evolution of the random graph $G(n,p)$ as $p
$ ``evolved'' from $0$ to $1$.  At $p=1/n$ a sudden and dramatic
change takes place in $G$.  When $p=c/n$ with $c<1$ the random $G$
consists of small components, the largest of size $\Theta (\log n)$.
But by $p=c/n$ with $c>1$ many of the components have ``congealed''
into a ``giant component'' of size $\Theta (n)$.  Erd\H{o}s and 
Renyi called this the {\em double jump}, the terms phase transition
(from the analogy to percolation) and Big Bang have also been proferred.
\par Now imagine an observer who can only see $G$ through a logical
fog.  He may refer to graph theoretic properties $A$ within a limited
logical language.  Will he be able to detect the double jump?  The
answer depends on the strength of the language.  Our rough answer
to this rough question is: the double jump is not detectible in the
First Order Theory of Graphs but it is detectible in the Second
Order Monadic Theory of Graphs.  These theories will be described
below.  We use the abbreviations {\em fotog} and {\em somtog} for
these two theories respectively. \par For any property $A$ and
any $c>0$ we define \[f(c)=f_A(c)=\lim_{n\rightarrow \infty} Pr[G(n,c/n)
\models A] \] Here $G\models A$ means that $G$ satisfies property $A$.
Beware, however, that we cannot presuppose the
existence of $f(c)$ as the limit might not exist.
\\ {\bf Theorem 1.}  Let $A$ be a {\em fotog} sentence.  The $f_A(c)$
exists for all $c$ and $f_A$ is an infinitely differentiable function.
Moreover $f_A$ belongs to the minimal family of functions $\cal F$
which contain the functions $0,1$ and $c$, are closed under addition,
subtraction, and multiplication by rationals and are closed under
base $e$ exponentiation so that $f\in { \cal F} \Rightarrow e^f \in
{\cal F}$.
\\ Examples. Let $A$ be ``there exists a triangle''. Then $f(c)=e^{-c^3/6}$.
Let $B$ be ``there exists an isolated triangle'', i.e., a triangle with
none of the three vertices adjacent to any other vertices besides
themselves.  Then\[f_B(c)=e^{-c^3e^{-3c}/6}\] Finally, call a triangle
$x,y,z$ unspiked if there is no point $w$ which is adjacent to
exactly one of $x,y,z$ and no other point.  Let $C$ be the property
that there is no unspiked triangle.  Then 
\[f_C(c)=e^{-c^3e^{-3ce^{-c}}/6}\]
\par Since $1$ is not a special value for these $f$ we say, roughly,
that the double jump is not detectible in fotog.  Our remaining 
results all concern somtog.
{\bf Theorem 2.} There is a somtog $A$ with
\[f_A(c)=\left\{ \begin{array}{ll} 0 & \mbox{if $c<1$}\\1 & \mbox{if $c>1$}
\end{array} \right. \]
\\{\bf Theorem 3.} For all somtog $A$ and $c<1$ the value $f_A(c)$ is
well defined. \\{\bf Theorem 4.} For all $c_0>1$ there is a somtog $A$
with $f_A(c)$ not defined for $c \geq c_0$.
\\{\bf Theorem 5.}For any $c<1$ and $\epsilon >0$ there is a decision
procedure that will determine $f_A(c)$ within $\epsilon$ for any
somtog $A$.
\\{\bf Theorem 6.} Let $c>1$.  Then there is no decision procedure that
separates the somtog $A$ with $f_A(c)=0$ from those with $f_A(c)=1$.
\par Certainly the situation with c=1 is most interesting but we do
not discuss it in this paper.
\\ {\em Description of Theories:} The First Order Theory of Graphs ({\em 
 fotog}) consists of an infinite number of variable symbols $(x,y,z\ldots
)$, equality $(x=y)$ and adjacency (denoted $x\sim y$) symbols, the usual
Boolean connectives ($\wedge ,\vee ,\neg \ldots$) and universal and
existential quantification ($\forall_x,\exists_y,\ldots$) over the
variables which represent vertices of a graph. Second Order Monadic
Theory of Graphs ({\em somtog}) also includes an infinite number of set
symbols($S,T,U\ldots$) which represent subsets of the vertices and
membership ($\in$) between vertices and sets ($x\in S$).  The set
symbols may be quantified over ($\forall_S,\exists_T\ldots$) as well
as the variables.  As an example, in fotog we may write
\[\forall_x\forall_y\exists_z\exists_w[x\sim z\wedge z\sim w\wedge w\sim
y]\] which means that all pairs of vertices are joined by a path of length
three.  However, one cannot say in fotog that the graph is connected. In
somtog we define $path(x,y,S)$ to be that $x,y\in S$ and every $z\in S$
is adjacent to precisely two other $w\in S$ except for $x$ and $y$ which
are each adjacent to precisely one point of $S$.  This has the
interpretation that $S$ gives an induced path between $x$ and $y$.  The
statement $\exists_Spath(x,y,S)$ holds if and only if $x$ and $y$ lie
in the same component since if they do a minimal path $S$ would be an
induced path.  We write $conn(x,y,R)$ for $\exists_SS\subset R\wedge
path(x,y,S)$ which means that in the restriction to $R$, $x$ and $y$
lie in the same component.  The property ``$G$ is connected'' is
represented by the somtog sentence
$\forall_x\forall_y\exists_Spath(x,y,S)$.
This ability to express $x$ and $y$ being joined by some path of arbitrary
size seems to give the essential strength of somtog over fotog. Now we
can prove Theorem 2.  Let $circ(S)$ be the sentence that $S$ is
connected and that every $v\in S$ is adjacent to precisely two $w\in S$.
Consider the sentence
\[A: \exists_{S,T,R}circ(S)\wedge circ(T)\wedge S\cap T=\emptyset
\wedge S\subset R\wedge T\subset R\wedge \forall_{x,y\in R}conn(x,y,R)\]
This has the interpretation that the graph contains a component ($R$)
with two disjoint circuits.  For this $A$ it is well known that
$f_A(c)=0$ when $c<1$ and $f_A(c)=1$ when $c>1.$ 
\section{The First Order World}
The results of this section were done independently and in more complete
detail in Lynch[1].  Here we attempt to give a more impressionistic 
picture of $G(n,c/n)$ through First Order glasses.
\par What does $G(n,p)$ look like?  To begin with, there
are lots of trees.  More precisely, for any tree $T$ and any $r$ there
are almost surely more than $r$ copies of $T$ as components of the graph.  
(This includes the trivial case where $T$ is a single vertex.)  What about
more complicated structures.  Let $B(x,R)$ denote the set of vertices
within distance $R$ of $x$, where we use shortest path as the metric.
The veracity of a fotog $A$ depends only on the values $B(x,R)$ for a 
fixed (dependent on $A$) $R$. (This is most certainly not the case
for somtog.)  For any fixed $c$ and $R$ a.a. all $B(x,R)$ will be
either trees or unicylic graphs.  (For $c>1, G(n,c/n)$ will have many
cycles in the giant component but they will be far apart.) To make
things a bit bigger let's define ${\cal H}={\cal H}_R$ to be the set
of graphs $H$ consisting of a cycle of size at most $R$ and trees of
depth at most $R$ rooted at each vertex of the cycle.
For $C\subset G$ let $B(C,R)$ denote the set of vertices
within distance $R$ of some $x\in C$.  For any 
$H\in {\cal H}$ let $X_H$ denote the number of cycles $C$
with $B(C,R)\cong H$.When $H$ has $v$ vertices, and hence
$v$ edges, with $w$ vertices not at depth $R$
\[ E[X_H]={n\choose v}p^v(1-p)^{wn}/\mid Aut(H)\mid  \sim \lambda_H\]
where
\[\lambda_H=c^ve^{-wc}/v! \mid Aut(H)\mid \]
Moreover, the $X_H$ act as independent Poisson distributions in
that for any $H_1,\ldots H_s\in {\cal H}$
\[lim_{n\rightarrow \infty} Pr[X_{H_i}=c_i, 1\leq i \leq s] =
\prod_{i=1}^s \lambda_{H_i}^{c_i}e^{-\lambda_{H_i}}/c_i!\]
\par The values $B(C,R)$ can be generated as follows.  For each
$3\leq i \leq R$ the number of $i$-cycles is Poisson with mean
$c^i/2i$.  We generate these cycles and then from each vertex
generate a pure birth process with Poisson $c$ births.
\par To check the veracity of a fotog $A$ one needs only examine
$B(x,R)$ and further one doesn't need to be able to count higher
than $R$.  Lets say two numbers are $R$-same if they are either 
equal or both at least $R$.  Say two rooted trees of depth $1$
are $R$-equivalent if the degrees of their roots are $R$-same.
Clearly there are $R+1$ equivalence classes. Suppose by induction
$R$-equivalent has been defined on rooted trees of depth $i$.
A rooted tree of depth $i+1$ may be identified naturally with
a set of rooted trees of depth $i$.  We call two rooted trees of
depth $i+1 R$-equivalent if every $R$-equivalence class of 
rooted trees of depth $i$ appears the $R$-same number of times
as a subtree.  We say $H,H'\in {\cal H}$ are $R$-equivalent if
their cycles are the same size and the vertices can be ordered
about the cycles so that the rooted trees emenating from the
corresponding vertices are $R$-equivalent.  We say graphs
$G,G'$ are $R$-equivalent if every equivalence class of 
$H\in{\cal H}$ appears the $R$-same number of times as a
$B(C,R)$.  
\par We use the following result about fotog.  For every fotog
$A$ there is an $R$ with the following property.  Let $G,G'$
be $R$-equivalent, both with all $B(x,R)$ either trees or
unicylic.  Suppose further that every $R$-equivalence class
of rooted trees of depth $R$ appears at least $R$ times as
a $B(x,R)$ in both graphs.  Then
\[G\models A \Longleftrightarrow G'\models A\]
\par With $R$ fixed and $c$ fixed the rooted trees above certainly
appear. By induction on the depth we may show that every $R$-equivalence
class of $H\in {\cal H}$ appears Poisson $\lambda (c)$ times in $G(n,c)$
where $\lambda \in {\cal F}$.  Each $R$-equivalence class of $G$ occurs
with probability $f(c)\in {\cal F}$. Then $Pr[A]$ is the finite sum 
of such $f(c)$ and is also in ${\cal F}$.
\section{The Second Order Monadic World Before
the Double Jump} Here we prove Theorem 5 and hence
the weaker Theorem 3.  Fix $c>0$ and $\epsilon >0$. There are $K,L$ so
that with probability at least $1-\epsilon$ all components of $G(n,c/n)$
are either trees or unicylic components of size at most $K$, and there
are less than $L$ of the latter.  Each possible family of unicyclic
components holds with a calculatable limit probability.  Knowing the
precise unicyclic components and that $G(n,c/n)$ has ``all'' trees
and no other components determines the veracity of $A$.  Hence $Pr[A]$
is determined within $\epsilon$.
\section{The Second Order Monadic World
After the Double Jump} We prove Theorem 6 by a reduction to the
Trakhtenbrot-Vought Theorem, which states that there is no decision
procedure which separates those fotog $A$ which hold for some finite
graph from those which do not. By a clean topological $T_k$($CTK_k$) in
a graph $G$ we mean an induced subgraph consisting of $k$ vertices,
one path between every pair of points, and nothing else. In \S 5 we
show for all $c>1$ and  all integers $k$ that $G(n,c/n)$ almost surely
contains a $CTK_k$.
For any fotog $A$ we define a somtog $A^+$ of the form
\[A^+: \exists_{S,T,U} clean(S,T)\wedge A^* \]
Here $clean(S,T)$ represents that $S$ is the vertex set of a $CTK_k$
on $T$. That is
\\(i) $S\subset T$
\\(ii) Every $x,y\in S$ have a unique $T_{x,y}\subset T$ with $part(x,y,T_{x,y})$
and $T_{x,y}\cap S=\{ x,y\} $.
\\(iii) There is no edge between any  $T_{x,y}$ and $T_{x',y'}$ except at the endpoints.
\par Now we transform $A$ to $A^*$ by 
\\(i) replacing $\forall_x$ and $\exists_y$ by $\forall_{x\in S}$ and 
$\exists_{y\in S}$
\\(ii) replacing $x\sim y$ by $T_{x,y}\cap U\neq \emptyset$, with $T_{x,y}$
defined  above.
\par If $A$ holds for no finite graph then $A^+$ holds for no finite graph.  
Suppose $A$ holds for a finite graph $H$  on, say, $k$ vertices $1\ldots k$. 
Almost surely $G$ contains a $CTK_k$ on vertices $T$ with endpoints $S$.  
Label $S$ by $x_1\ldots x_k$ arbitrarily.  Let $U$ consist of one vertex from
$T_{x_i,x_j}$ (not an endpoint) for each pair $\{ x_i,x_j\}$ with
$\{ i,j\} \in H$.  Then $A^*$ holds.  That is, $A^+$ holds a.a.
\par A decision procedure that could separate somtog $B$ with $f_B(c)=1$ from
those with $f_B(c)=0$ could, when applied to $B=A^+$, be used to determine
if $A$ held for some finite graph, and this would contradict the 
Trakhtenbrot-Vought Theorem.
\\ {\em Nonconvergence.} To prove Theorem 4 we use a somewhat complicated
graph.  Let $k_1$ be a positive real and $K$ a positive integer.  ($k_1=5,
K=100$ is a good example.)  We define, for all sufficiently large $n$, a
graph $H=H(k_1,K,n)$.  Let $w$ be the nearest integer to $k_1\log n$
divisible by $K$ (a technical convenience) so that $w\sim k_1\log n$.
(Asymptotics are in $n$ for fixed $k_1,K$.  Begin with two points $S0,S1$
and three vertex disjoint paths, each of length $w$, between them.  Call
this graph $H^-$.  Let $AR$ (which stands for arithmetizable) consist of
every $K$-th vertex on each of the paths, excluding the endpoints. Thus
$AR$ will have $\frac{w}{K}-1$ points from each path, a total of
$w_1=3[\frac{w}{K}-1]$ points.  Order the three paths arbitrarily and
order the points of $AR$ on a path from $S0$ to $S1$ so that the points
of $AR$ are labelled $1,\ldots,w_1$.  Now, using this labelling,
between every pair $i,2i$ add a path of length $w$.  (These paths all
use new vertices with no additional adjacencies.)  Now between every pair
$i,2^i$ add a path of length $w$. Now between every pair $i,tower(i)$
add a path of length $w$.  (The function $tower(i)$ is defined inductively
by $tower(1)=2$,$tower(i+1)=2^{tower(i)}$.)  Finally between every pair
$i,wow(i)$ add a path of length $w$.  (The function $wow(i)$ is defined
inductively by $wow(i)=2$, $wow(i+1)=tower(wow(i))$.)  This completes
the description of the graph $H=H(k_1,K,n)$.  
\par In \S 5 we prove that for every $c>1$ there exist $k_1,K$ so that
$G(n,c/n)$ almost surely contains an induced copy of $H=H(k_1,K,n)$.
We assume that here, and with $H$ in mind construct a somtog 
sentence $A=A_K$ which shows nonconvergence.
\par The sentence $A=A_K$ will be built up in stages.  First we say
there exist vertices $S0,S1$ and sets $P_1,P_2,P_3$ so that each
$P_i$ gives a path from $S0$ to $S1$, the $P_i$ overlap only at $S0,S1$,
and there are no edges between $P_i$ and $P_j$ except at the endpoints.
Second we say there exists a set $AR\subset P_1\cup P_2\cup P_3-\{
S0,S1\}$ so that for any path $x_1\cdots x_K$ in any $P_i$ that
$AR$ contains exactly one of the $x_1,\ldots,x_K$.  (Here the sentence
depends on the choice of the fixed integer $K$.)  We define an auxilliary
binary relation $<$ on $AR$.  If $x,y\in AR\cap P_i$ we define $x<y$ by
the existence of a subset of $P_i$ which is a path from $S0$ to $x$ which
does not contain $y$.  If $x\in AR\cap P_i$ and $y\in AR\cap P_j$ with
$i\neq j$ we define $x<y$ to be $i=1,j=2$ or $i=1,j=3$ or $i=2,j=3$. On
$AR$ we define the auxilliary binary predicate $next(i,j)$ by $i<j$ and
there does not exist $k\in AR$ with $i<k$ and $k<j$. We define the 
unary predicate $ONE(i)$ by $i\in AR$ and there is no $j<i$ and $TWO(i)$
by $i\in AR$ and $j<i\leftrightarrow ONE(j)$.  We say there are unique
$i,j$ with $ONE(i)$, $TWO(j)$.  For convenience we write $1,2$ for these
elements henceforth.
\par Now to arithmetize $AR$.  We say there exists vertex sets $DOUBLE$,
$EXP$, $TOWER$ and $WOW$.  We define auxilliary binary predicate $double$
on $AR$ by $double(x,y)$ if $x<y$ and there is a path from $x$ to $y$
in $DOUBLE$; and we similarly define binary predicates $exp$, $tower$
and $wow$.  We say $double(1,2)$ and $double(x,y)\cap double(x,z)\rightarrow
y=z$ and $double(x,y)\cap next(x,x_1)\cap next(y,y_1)\cap next(y_1,y_2)
\rightarrow double(x_1,y_2)$ and if $double(x,y)$ and $next(x,x_1)$ and
there do not exist
$y_1,y_2$ with $next(y,y_1)\cap next(y_1,y_2)$ then there does not exist
$z$ with $double(x_1,z)$ and finally if $double(x,y)$ and $x'<x$ then
there exists $y'$ with $double(x',y')$. We say $exp(1,2)$ and $exp(x,y)\cap
exp(x,z)\rightarrow y=z$ and $exp(x,y)\cap next(x,x_1)\cap double(y,y_1)
\rightarrow exp(x_1,y_1)$ and if $exp(x,y)$ and $next(x,x_1)$ and there
does not exist $y_1$ with $double(y,y_1)$ then there does not exist $z$
with $exp(x_1,z)$ and finally if $exp(x,y)$ and $x'<x$ then there exists 
$y'$ with $exp(x',y')$.  The properties for $tower$ are in terms of $exp$
exactly as the properties for $exp$ were in terms of $double$ and the
properties for $wow$ are in terms of $tower$ in the same way.
\par On $AR$ we define unary predicates $even(x)$ by $\exists_y double(y,x)$
and $invwow(x)$ by there existing $y$ with $wow(x,y)$ but for all $x'>x$ there
do not exist $y'$ with $wow(x',y')$.  The sentence $A=A_K$ concludes
by saying there exists $x$ with $even(x)\cap invwow(x)$.
\par Now we show that $\lim Pr[G(n,p)\models A_K]$ does not exist, moreover
that the $\limsup$ is one and the $\liminf$ is zero. On the integers define
$wow^{-1}(y)$ to be the biggest integer $x$ with $wow(x)\leq y$.  First
let $n\rightarrow\infty$ through that subsequence for which $wow^{-1}(w_1)$
is even. (Recall $w_1=\Theta(\log n)$ was the size of $AR$.)  Suppose $G(n,p)$
contains an induced copy of $H$. ($k_1,K$ depend only on $c$ and so are
already fixed.)  On $H$ there do exist the vertices $S0,S1$,
the sets $P_1,P_2,P_3$, $AR$, $DOUBLE$, $EXP$, $TOWER$ with all the properties
of $A_K$.  (Indeed, $A_K$ was created with that in mind.)  Under the labelling
$1,\ldots, w_1$ the predicates $double,\ldots$ correspond to the actual
numbertheoretic predicates and the $x=wow^{-1}(w_1)$ has $invwow(x)$ and
$even(x)$ so $A_K$ holds.  But $G(n,p)$ contains an induced copy of $H$
almost surely so the limiting probability on this subsequence is one.
\par In the other direction, let $n$ go to infinity through a subsequence
with the property that for all $m$ with (leaving some room) $\log\log n<m<n$
the value $wow^{-1}(m)$ is odd.  (Such $n$ exist since $wow{-1}$ is constant
for such a long time.)  Here is the crucial random graph fact:  There is
a $\delta=\delta(c)$ so that in $G(n,c/n)$ almost surely all subconfigurations
consisting of two vertices and three paths between them have size at least
$\delta\log n$.  (This uses a simple expectation argument.  The number of
configurations of $t$ vertices and $t+1$ edges giving the above graph is
$O(n^tp^{t+1})=O(c^{t+1}/n)=o(1)$ when $t<\delta\log n$.)  Thus almost
surely any $AR$ that satisfies the conditions of $A_K$ will have
$|AR|=m>\delta'\log n>\log\log n$.  The conditions on $double,\ldots$
{\em force} $AR$ to be arithmeticized so that $\exists_x invwow(x)\cap
even(x)$ will not occur when $wow^{-1}(m)$ is odd.  Thus almost surely
$A_K$ will not be satisfied.
\section{A Variance Calculation}
We fix $c>1$, set $p=c/n$ and let $G\sim G(n,p)$.  We consider a graph
$H=H(k_1,K,n)$ as defined in \S 4. We give a description of $H$ suitable for
our purposes.  Set $w=\lceil k_1\log n \rceil$.  Take two vertices and
draw three vertex disjoint paths each of length $w$.  This gives a 
graph $H^-$.  On $H^-$ a set of  pairs of vertices
$\{a,a'\}$ are specified, no $a$ lying in more than eight such pairs.
We let $l$ denote the precise number of such pairs so that $l\sim
\epsilon\log n$.  By making $K$ large we can make $\epsilon$ as small
as desired.
Between each such pair a path of length $w$ is placed with new vertices.
This gives the graph $H$.  It has $v=\Theta(\ln^2n)$ vertices and $e$
edges where $e=3w+lw\sim \epsilon k_1\log^2n$ and
$e-v=l+1\sim\epsilon\log n$. Let us denote the vertices of $H$
by $1,\ldots,v$.
\subsection{The Second Moment Method}
Our object in this section is to show that, for appropriate
$k_1,\epsilon$, the random $G(n,p)$ almost surely contains an induced
copy of $H$. Let $X$ be the number of $v$-tuples $(a_1,\ldots,a_v)$
of distinct vertices of $G$ so when $\{i,j\}\in E(H)$ then $\{a_i,a_j\}
\in E(H)$.  That is, $X$ is a count of copies of $H$ in $G$ though these
copies may have extra edges and a given copy may be multiply counted if
$H$ has automorphisms. Clearly
\[ E[X]=(n)_vp^e\sim n^vp^e=c^e/n^{e-v}  \]
which is
\[ n^{(\epsilon\log n)(k_1\log c-1+o(1))} \]
from the estimates above.  We first require that
\[ k_1\log c > 1  \]
which assures that $E[X]$ is a positive power of $n^{\log n}$.  The
crucial calculation will be to show
\[ Var[X]=o(E[X]^2) \hspace*{2cm} (V1)  \]
>From this, by Chebyschev's Inequality $X>.99 E[X]$ (say) almost surely. True,
$X$ counts noninduced copies of $H$.  But let $X^+$ be a count of
all copies of any $H^+$ consisting of $H$ with one additional edge added.
There are $\Theta(\log^4n)$ choices of that edge and for a given choice
the expected number of such copies is $pE[X]$ so that $E[X^+]=
O(\log^4n/n)E[X]=o(E[X])$ and so by Markov's Inequality almost surely
$X^+<E[X]/2$, say.  So almost surely there are more than $.99E[X]$
copies of $H$ and fewer than $.5E[X]$ total copies of graphs containing
$H$ and one more edge so therefore there is at least one copy of $H$ with
no additional edge, i.e., the desired induced copy.
\par Hence it suffices to show (V1).
\par {\em Remark.}  To illustrate the complexities suppose we condition
$G(n,p)$ on a fixed copy of $H^-$ and let $Z$ be the expected number of
extensions to $H$.  The expectation argument above gives that 
$E[Z]\sim (n^{w-1}p^w)^l = (c^w/n)^l$ which is $n^{\Theta(\log n)}$.
However for there to be any extensions each of the at least $l/8$
vertices of $H^-$ that is supposed to have a path coming out of
it must have at least one edge besides those of $H^-$.  Any particular
vertex fails this condition with probability $e^{-c}$ and these events
are independent so that the probability that $Z\neq 0$ is bounded from
above by $(1-e^{-c})^{l/8}$ which is polynomially small.  This illustrates
that the expected number of thingees being large does not
{\em a priori} guarantee that almost surely there is a thingee.
\par Of course, (V1) is equivalent to showing $E[X^2]\sim E[X]^2$.
By the symmetry of copies, $E[X^2]$ is $E[X]$ times the expected
number of copies of $H$ given the existence of a particular copy 
of $H$.  Let
\[ V=\{1,\ldots,m\} \]
and let us specify a particular copy of $H$ on vertex set $V$
with $1,\ldots,3w-1$ being the vertices of $H^-$.
Let $G^*=G^*(n,p)$ be the random graph on vertex set
$1,\ldots,n$ where for $i,j\in V$ and $\{i,j\}\in E(H)$
we specify that $\{i,j\}\in E(G)$ but all other pairs $i,j$ are
adjacent in $G$ with independent probabilities $p$.  (Note that
even those $i,j$ with $i,j\in V$ but $\{i,j\}\not\in E(H)$
have probability $p$ of being in $G^*$.)  Let $E^*[X]$ denote
the expectation of $X$ in $G^*$.  Then it suffices to show
\[ E^*[X]\sim E[X] \hspace*{2cm} (V2)  \]
\subsection{The Core Calculation: Expectation for Paths}
We shall work up to $E*[X]$ in stages.  Let $P_s(a,b)$ denote
the expected number of paths of length $s$ between vertices $a,b$,
with the graph distribution $G^*(n,p)$.  (As a benchmark note that
in $G(n,p)$ this expectation would be $(n-2)_{s-1}p^s\sim n^{s-1}p^s$.)
$P_s(a,b)$ is simply the sum over all tuples $(a_0,\ldots,a_s)$ with
$a_0=a,a_s=b$ of distinct vertices of $G$ of $p^{\alpha}$ where 
$\alpha$ is the number of edges of the path $a_0\cdots a_s$ which are
{\em not} in $H$.  Let $P_s^-(a,b)$ denote the expected number of
such paths where we further require that $a$ is not adjacent to $a_1$
in $H$.  (When $a\not\in V$ these are the same.)  We shall define
inductively $x_s,x_s^-$ which provide upper bounds to $P_s(a,b)$ and
$P_s^-(a,b)$ respectively under the further assumption that $b\not\in V$.
(We shall see that $P_s^-(a,b)$ is dominated by paths which do not
overlap $H$ but that for $P_s(a,b)$  there is a contribution from
those paths which are paths in $H$ for their initial segment.)
Clearly we may set $x_1=x_1^-=p$. Let $x_s,x_s^-$ satisfy the following:
\[ x_s^- = pnx^-_{s-1}+pmx_{s-1}  \]
\[ x_s = x_s^-+\sum_{k=1}^{s-1}50kx^-_{s-k} \]
We claim such $x_s,x_s^-$ provide the desired upper bounds.  To bound
$P_s^-(a,b)$ split paths $aa_1\cdots a_{s-1}b$ according to $a_1\in V$
($\leq m$ possibilities) and $a_1\not\in V$ ($\leq n$ possibilities).
Note we are excluding the case where $a,a_1$ are adjacent in $H$.
For a given $a_1$ the expected number of paths is $pP_{s-1}(a_1,b)$
(as we must have $a,a_1$ adjacent). When $a_1\not\in V$ this is by
induction at most $px^-_{s-1}$ and when $a_1\in V$ this is by 
induction at most $px_{s-1}$ so $P_s^-(a,b)\leq x_s^-$ by induction.
Bounding $P_s(a,b)$ is a bit more complex.  Those paths for which
$a,a_1$ are not adjacent in $H$ contribute at most $x_s^-$ by
induction.  Otherwise, let $k$ be the least integer for which 
$a_k,a_{k+1}$ are not adjacent in $H$. (As $b\not\in V$ this is
well defined and $1\leq k<s$.) We pause for a technical calculation.
\par We claim that in $H$ for any $k\leq w$ there are at most $50k$
paths of length $k$ beginning at any particular vertex $v$. 
Suppose $a\in H^-$. There are at most four such paths staying in $H^-$.
Once leaving $H^-$ the path is determined (since critically $k\leq w$,
the path length) and there are at most $8k$ ways of determining when
and how to leave $H^-$.  The argument with $a\not\in H^-$ is similar,
we omit the details.  Of course $50k$ is a gross overestimate but
we only use that it is a $O(k)$ bound.
\par Back to  bounding $P_s(a,b)$. For a given $k$ there are at most
$50k$ choices for $a_1\cdots a_k$ and fixing those there is a 
contribution of $P_{s-k}^-(a_k,b)\leq x_{s-k}^-$ to $P_s(a,b)$.  Thus
$P_s(a,b)\leq x_s$ by induction.
\par Now to bound the values $x_s,x_s^-$ given by the inductive formulae.
Let $L$ be fixed (dependent only on $c$) so that
\[ L > 1 + \sum_{k=1}^{\infty} 50kc^{-k} \]
and set
\[ X_s = Lp^sn^{s-1} \]
\[ X_s^-=p^sn^{s-1}(1+L\frac{ms}{n}) \]
We claim $x_s\leq X_s$ and $x_s^-\leq X_s^-$.  For this we merely check
(recall $pn=c$)
\[ X_s^-+\sum_{k=1}^{s-1}50kX_{s-k}^- \leq (1+L\frac{ms}{n})p^sn^{s-1}
(1+\sum_{k=1}^{s-1}50kc^{-k}) < Lp^sn^{s-1}=X_s  \]
as $Lms/n=o(1)$ and that
\[ pnX_{s-1}^-+pmX_{s-1}=X_s^- \]
Thus we have shown
\[ P_w(a,b)\leq Ln^{w-1}p^w  \]
when $b\not\in V$ and further
\[ P_w(a,b) \leq n^{w-1}p^w[1+O(\frac{\log^3n }{n})]  \]
when $a,b\not\in V$.  
\par Now (thinking of $a,b\in V$) we seek a general bound $y_s$ for
$P_s(a,b)$. We set $y_1=1$ (as perhaps $a,b$ are adjacent in $H$)
and define inductively
\[ y_s=3+pnx_{s-1}+pmy_{s-1}+50sp + \sum_{k=1}^{s-2}50k[pnx_{s-k-1}+
pmy_{s-k-1}] \]
We claim $P_s(a,b)\leq y_s$   for $1\leq s\leq w$.  Of the potential
paths $aa_1\cdots a_{s-1}b$ there are at most three which are paths in
$H$ and they contributes at most three.  There are less than $50s$ cases
where $aa_1\cdots a_{s-1}$ is a path in $H$ but $a_{s-1},b$ are not 
adjacent in $H$ and they each contribute $p$.
The cases with $a_1\not\in V$
contribute at most $pnx_{s-1}$. The cases with $a_1\in V$ but not
adjacent to $a$ in $H$ contribute at most $pmy_{s-1}$.  Otherwise
let $1\leq k\leq s-2$ be the least $k$ so that $a_k,a_{k+1}$ are
not adjacent in $H$.  There are at most $50k$ choices of $a_1\cdots a_k$.
Then there are at most $n$ choices of $a_{k+1}\not\in V$ and each
contributes $px_{s-k-1}$ and at most $m$ choices of $a_{k+1}\in V$
and each contributes $py_{s-k-1}$.  
\par Now fix a constant $M$ satisfying
\[ M > M_1=L[1+\sum_{k=1}^{\infty}50kc^{-k}]  \]
We claim that for $1\leq s\leq w$
\[ P_s(a,b) \leq 4 + Mp^sn^{s-1}  \]
By the previous bounds on $x_s$ we bound
\[ pnx_{s-1}+\sum_{k=1}^{s-2}50kpnx_{s-k-1} < M_1p^sn^{s-1} \]
We bound $3+50sp<3.01$. By induction we bound
\[ pmy_{s-1}+\sum_{k=1}^{s-2}50kpmy_{s-k-1} < 50s^2pm[4+Mp^sn^{s-1}]
< .01+(M-M_1)p^sn^{s-1} \]
since $50s^2pm=O(\log^4n/n)=o(1)$, completing the claim.  We are really
interested in the case $s=w$.  Note $p^wn^{w-1}=c^w/n$ is asymptotically
a positive power of $n$ by the choice made of $k_1$ earlier. Thus
the $+4$ may be absorbed in $M$ and we have that
\[ P_w(a,b) < Mp^wn^{w-1}  \]
for all $a,b$ while if  $a,b\not\in V$ then we have
the better bound
\[ P_w(a,b) < p^wn^{w-1}[1+O(\frac{\ln^3n}{n})]  \]
\subsection{Expectation of Copies of $H$}
Now we turn to the full problem of bounding $E^*[X]$. 
Recall we have labelled  $H$
so that $1,\ldots,3w-1$ are the vertices of $H^-$.  Recall  $l$ 
denotes the number of $w$-paths in going from $H^-$ to $H$ and
recall $l\sim \epsilon\log n$.
$E^*[X]$ is the
sum over all $m$-tuples $(a_1,\ldots,a_m)$ of distinct vertices of
the probability (in $G^*(n,p)$) that these $a$s (in this order) give 
a copy of $H$.  For each $a_1,\ldots,a_{3w-1}$ the contribution
of $m$-tuples with this start is bounded from above by
\[ p^{\alpha}\left[Mp^wn^{w-1}\right]^{l-A}\left[p^wn^{w-1}[1+
O(\frac{\ln^3n}{n})]\right]^A  \]
Here $\alpha$ is the number of adjacencies $i,j$ in $H^-$ with
$a_i,a_j$ not adjacent in $H$.  $A$ is the number of pairs $i,j$
in $H^-$ which are joined in $H$ by a $w$-path and for which neither
$a_i$ nor $a_j$ is in $V$.  $l-A$ is then the  remaining number of pairs $i,j$
in $H^-$ joined in $H$ by a $w$-path.  
\par To see this note that for fixed $a_1,\ldots,a_{3w-1}$ and any choice
of $w$-paths $P_1,\ldots,P_l$ that are vertex disjoint the probability that
they are all paths in $G$ is simply the product of the probabilities for each
path.  Adding over all $P_1,\ldots,P_l$ is then at most the product over $j$
of adding the probabilities for each $P_j$, and these are precisely  what the
bracketed terms bound.  The $p^{\alpha}$, of course, is the probability that
the $a_i$ have the proper edges of $H^-$.
\par Now we split the contribution to $E[X^*]$ into two classes.  First
consider all those $a_1,\ldots,a_{3w-1}$ with {\em no} $a_i\in V$.  There
are at most $n^{3w-1}$ such tuples and each gives $p^{3w}$ with
$3w$ being the number of edges in $H^-$.  For each $A=l$ so this gives a
\[ \left[p^wn^{w-1}[1+O(\frac{\ln^3n}{n})]\right]^l  \]
factor.  As $l=O(\ln n)$
\[ [1+O(\frac{\ln^3n}{n})]^l = 1+o(1)  \]
so this entire contribution is asymptotic to $n^{3w-1+l(w-1)}p^{3w+lw}$
which is asymptotic to $E[X]$, the expectation in $G(n,p)$.  That is,
the main contribution (among the $a_1\cdots a_{3w-1}$ that don't overlap
$V$) to $E^*[X]$ is by those copies of $H$ that don't overlap $H$ at all.
To show (V2) it now suffices to show that the remaining contributions
to $E^*[X]$ are $o(E[X])$.
\par There are $O(\ln^3n)$ choices of a pair $i\leq 3w-1$ and $a_i\in V$.
Fix such a pair and consider the contribution to $E^*[X]$ with $a_i$
this fixed value.  In $H^-$ fix a cycle $C$ of length $2w$ that $i$ lies
on. The expected number of cycles of length $2w$ in $G$ through $a_i$
bounded from above by $P_{2w}(a_i,a_i)$ which we've shown is at most
$Mn^{2w-1}p^{2w}$.  (The analysis done for $P_s$ for $1\leq s\leq w$
extends with no change to $s=2w$.)  Let $Q,Q'$ be the points
of $H^-$ of degree three. Fixing such a cycle  the points $a_Q,a_{Q'}$
are now fixed and the expected number of paths of length $w$ between
them is at most $P_w(a_Q,a_{Q'})\leq Mn^{w-1}p^w$.  Together the 
expected number of extensions of $a_i$ to a copy of $H^-$ is less than
$O(n^{3w-2}p^{3w})$, off from the expected number of copies of $H^-$
by a factor of $n^{-1+o(1)}$  The $O(\ln^3n)$ factor of the choices
of $i,a_i$ can be absorbed into to $o(1)$ so that the expected number
of copies of $H^-$ overlapping $V$ is still only $n^{-1+o(1)}$ times
the expected number of copies of $H^-$ in $G(n,p)$. Now given
a copy $a_1,\ldots,a_{3w-1}$ of $H^-$ the expected number of extensions
to $H$ in $G^*(n,p)$ is at most $M^l$ times what it is in $G(n,p)$,
the extreme case when all  $A=0$, e.g., all $a_i\in V$.  Thus the
total contribution to $E^*[X]$ from copies in which $H^-$ overlaps
$V$ is at most 
\[ n^{-1+o(1)}M^lE[X]  \]
Recall $l\sim\epsilon\log n$.  Up to now all constants $k_1,L,M$ have
depended only on $c$ and not $\epsilon$.  Now (and formally this is at
the very start of the proof, in the definition of $H$) we fix $K$ so 
large that $\epsilon$ is so small so that
\[ \epsilon(\log M) < 1  \]
This assures that $n^{-1+o(1)}M^l$ is $n$ to a negative power.  Thus
this contribution to $E^*[X]$ is only $o(E[X])$.  Hence 
$E^*[X]\sim E[X]$ which concludes the argument.         
\subsection{Clean Topological $k$-Cliques}
For the proof of Theorem 6 we require that for every $c>1$ and every
integer $k$ that $G(n,c/n)$ almost surely contains a $CTK_k$.  Fix
$c,k$.  We fix a real $k_1$ with
\[ k_1\log c > 1 \]
Set $w=\lceil k_1\ln n\rceil$.  We define $H=H(k,k_1,n)$ to consist
of $k$ ``special'' points and between each pair of special points
a path of length $w$.  Set $v=k+{k\choose 2}(w-1)$, the number of
vertices and $t={k\choose 2}-k$, $e=v+t$ so that $e$ is the number
of edges.  Note $e,v=(k_1{k\choose 2}+o(1))\log n$.
We show that almost surely $G(n,c/n)$ contains a copy
of $H$.  As the argument is very simpler (and simpler) than that
just given, we shall give the argument in outline form.  Letting
$X$ denote the number of copies of $H$ we have
\[ E[X]=(n)_vp^e\sim n^vp^e=c^en^{-t}=n^{k_1{k\choose 2}\log c-t+o(1)}
\]
which is a positive power of $n$.  Now we need show $E^*[X]\sim
E[X]$ where $E^*[X]$ is the expected number of copies of $H$ conditioning
on a fixed copy of $H$.  Let us specify the fixed copy to be on vertex set
$V=\{1,\ldots,v\}$ with $1,\ldots,k$ being the special vertices and 
let $G^*(n,p)$ be $G(n,p)$ conditioned on this copy.  As before
we let $P_w(a,b)$ be the expected number of paths of length $w$ between
$a,b$ in $G^*$. Then, as before, there is a constant $M$ so that
\[ P_w(a,b) < Mn^{w-1}p^w  \]
for all $a,b$ while
\[ P_w(a,b) = n^{w-1}p^w(1+o(1))  \]
if either $a$ or $b$ is not in $V$. (We will not need the more precise
error bound for this problem.)
\par We split the contribution to $E^*[X]$ into two groups. The 
$a_1,\ldots,a_k$ which do not overlap $V$ contribute
\[ n^k\left[ n^{w-1}p^w(1+o(1))\right]^{{k\choose 2}} \]
to $E^*[X]$.  Since $k$ is {\em fixed} this is asymptotically $n^vp^e$
which is asymptotically $E[X]$.  There are only $n^{k-1+o(1)}$
different $a_1,\ldots,a_k$ which do overlap $V$.  For each the
contribution to $E^*[X]$ is at most
\[ \left[Mn^{w-1}p^e\right]^{{k\choose 2}}  \]
Since $M$ and $k$ are constants this only a constant times the 
contribution to $E[X]$.  Thus the total contribution from these
intersecting $a_1,\ldots,a_k$ is $n^{-1+o(1)}E[X]$ and thus
$E^*[X]\sim E[X]$ as required.

\vspace{1cm}

{\bf References.}
\\ 1. J. Lynch, Probabilities of Sentences about Very Sparse Random Graphs,
Random Structures and Algorithms (to appear)

\end{document}